\documentclass[10pt, a4paper]{article}
\usepackage{rotate}
\usepackage{url}
\author{Andrew Aberdein\\
University of Edinburgh\\
{\tt andrew.aberdein@dunelm.org.uk}}
\title{The Informal Logic of Mathematical Proof}
\newif\ifpdf
\ifx\pdfoutput\undefined
\pdffalse % we are not running PDFLaTeX
\else
\pdfoutput=1 % we are running PDFLaTeX
\pdftrue
\fi

\ifpdf
\usepackage[pdftex]{graphicx}
\else
\usepackage{graphicx}
\fi

\ifpdf
\DeclareGraphicsExtensions{.pdf, .jpg}
\else
\DeclareGraphicsExtensions{.eps, .jpg}
\fi
\begin{document}
\maketitle

%\section{Applying informal logic to mathematics}
Paul Erd\H{o}s famously remarked that `a mathematician is a machine for turning coffee into theorems' \cite[p.~7]{Hoffman}.  The proof of mathematical theorems is central to mathematical practice and to much recent debate about the nature of mathematics. This paper is an attempt to introduce a new perspective on the argumentation characteristic of mathematical proof.  I shall argue that this account, an application of  informal logic to mathematics, helps to clarify and resolve several important philosophical difficulties.

It might be objected that formal, deductive logic tells us everything we need to know about mathematical argumentation.  I shall leave it to others \cite[for example]{Rav} to address this concern in detail.  However, even the protagonists of explicit reductionist programmes---such as logicists in the philosophy of mathematics and the formal theorem proving community in computer science---would readily concede that their work is not an attempt to capture actual mathematical practice.
Having said that, mathematical argumentation is certainly not inductive either.  Mathematical proofs do not involve inference from particular observations to general laws. 
A satisfactory account of mathematical argumentation must include deductive inference, even if it is not exhausted by it.  It must be complementary, rather than hostile, to formal logic.  My contention is that a suitable candidate has already been developed independently: informal logic.

Informal logic is concerned with all aspects of inference, including those which cannot be captured by logical form.  It is an ancient subject, but has been a degenerating research programme for a long time.  Since the nineteenth century it has been overshadowed by the growth of formal logic.  More fundamentally, it has suffered by identification with  the simplistic enumeration of fallacies, without any indication of the circumstances in which they are illegitimate.  Since most fallacies can be exemplified in some contexts by persuasive, indeed valid, arguments, this approach is of limited use.  In recent decades more interesting theories have been developed.  I shall look at two of the most influential, and discuss their usefulness for the analysis of mathematical proof.

\section{Toulmin's pattern of argument}
\begin{figure}[htb]
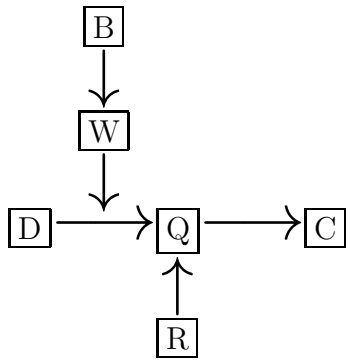

\caption{Toulmin's DWC pattern}
\label{fig:Toulmin}
%\unitlength 1in
%\begin{picture}(9.25,3.4)(0.8,0.3)
%\put(0,0){\includegraphics[angle=90, width=9.25in, height=4.2in]{ToulminPDF}}
%\end{picture}
\Huge
\[
\stackrel{\stackrel{\fbox{\large D} \stackrel{\stackrel{\stackrel{\stackrel{\fbox{\large B}}{\mbox{$\downarrow$}}}{\fbox{\large W}}}{\mbox{$\downarrow$}}}{\mbox{$\longrightarrow$}} \fbox{\large Q} \mbox{$\longrightarrow$} \fbox{\large C}
}{\mbox{$\uparrow$}}}{\fbox{\large R}} 
\]
\normalsize
\end{figure}

One of the first modern accounts of argumentation is that developed in Stephen Toulmin's \textit{The Uses of Argument} \cite{Toulmin58}.  Toulmin offers a general account of the layout of an argument, as a claim (C) derived from data (D), in respect of a warrant (W).  Warrants are general hypothetical statements of the form `Given D, one may take it that C' \cite[p.~99]{Toulmin58}.  Hence the laws of logic provide a warrant for deductive inferences.  However, the pattern is intended to be more general, and provides for different, weaker warrants, although these would not permit us to ascribe the same degree of certainty to C.  This is recognized by the inclusion of a modal qualifier (Q), such as `necessarily', `probably', `presumably',\dots, in the pattern.  If the warrant is defeasible, we may also specify the conditions (R) under which it may be rebutted.  Finally, the argument may turn on the backing (B) which can be provided for W.  Toulmin's claim is that the general structure of a disparate variety of arguments may be represented as in Figure~\ref{fig:Toulmin}.

Interpreting the letters as above, this diagram may be read as follows: ``Given D, we can (modulo Q) claim C, since W (on account of B), unless R''.  In Toulmin's vintage example: ``Given that Harry was born in Bermuda, we can presumably claim that he is British, since anyone born in Bermuda will generally be British (on account of various statutes \dots), unless he's a naturalized American, or his parents were aliens, or \dots''.  In simpler examples B, Q and R may not all be present, but D, W and C are taken to be essential to any argument, hence the description of this model as the DWC pattern.  Toulmin stresses the field dependency of the canons of good argument: what counts as convincing may vary substantially between the law court, the laboratory and the debating chamber.  In particular, what counts as acceptable backing will turn significantly on the field in which the argument is conducted \cite[p.~104]{Toulmin58}.

Toulmin's work has been very influential in the study of argument, despite an initially chilly reception amongst philosophers and logicians.\footnote{`Unanimous \dots condemnation' according to van Eemeren \textit{\& al}.\ \cite[p.~164]{vanEemeren+}, who survey the book's reviewers.}  It was quickly adopted by communication theorists, after the publication in 1960 of a celebrated paper by Wayne Brockriede and Douglas Ehninger \cite{Brockriede}.  Toulmin's account of argumentation is now the dominant model in this field.  More recently, his work has been widely studied by computer scientists attempting to model natural argumentation \cite[for example]{Newman}.  One purpose of Toulmin's critique of logic is to dispute the utility of formal logic for the analysis of any argumentational discourse \textit{other} than mathematics: he regards mathematical proof as one of the few success stories for the formal logic tradition.  Nevertheless, mathematical proofs can be subsumed under the DWC pattern, where the warrant is backed by various axioms, rules of inference and mathematical techniques providing grounds for supposing the claim to be necessary, given the data.  Toulmin provides an example in a later collaborative work \cite[p.~89]{Toulmin+}, by reconstructing Theaetetus's proof that there are exactly five polyhedra.  The data and warrant consist of various facts about the platonic solids, the warrant is backed by the axioms, postulates and definitions of three-dimensional Euclidean geometry, and the modal qualifier `with strict geometrical necessity' admits of no rebuttal or exception within the bounds of Euclidean geometry.

\section{Applying Toulmin to mathematics}

The significance of Toulmin's work for mathematical proof is explored at greater length in what I believe to be the only study so far of the application of informal logic to mathematics, a paper written in Catalan by Jes\'{u}s Alcolea Banegas \cite{Alcolea}.\footnote{I am grateful to Miguel Gimenez of the University of Edinburgh for translating this paper.}  Alcolea makes use of a further distinction of Toulmin's, introduced in \cite{Toulmin+}: that between \textit{regular} and \textit{critical} arguments.  This distinction echoes Kuhn's contrast between normal and revolutionary science: a regular argument is an argument within a field which appeals to the already well-established warrants characteristic of the field, whereas a critical argument is an argument used to challenge prevailing ideas, focusing attention on the assumptions which provide a backing for the warrants of regular arguments.  Critical arguments must therefore appeal to different warrants.  Mathematical proofs are regular arguments, although they may give rise to critical arguments if they are especially interesting or controversial.  Conversely, metamathematical debates are critical arguments, but they often provide new opportunities for proofs, that is, regular arguments.

Alcolea uses Toulmin's layout to reconstruct one regular and one critical argument from mathematics.  The critical argument, the debate over the admissibility of the axiom of choice, is the more fully developed and persuasive of Alcolea's case studies.  It is perhaps not too surprising that critical arguments in mathematics are similar to critical arguments in the other sciences, since ultimately they are not arguments \textit{in} mathematics, but arguments \textit{about} mathematics, that is to say they are metamathematical.  However, my concern is primarily with the argumentation of mathematics itself, rather than that of metamathematics.  Hence I shall concentrate on Alcolea's example of a regular argument: Kenneth Appel and Wolfgang Haken's proof of the four colour conjecture.\footnote{For further detail of the proof see \cite{Appel+}, \cite{Wilson} or \cite{MacKenzie}.}  He reconstructs the central argument of the proof as a derivation from the data $D_{1}$--$D_{3}$
\begin{quotation}
\noindent ($D_{1}$) Any planar map can be coloured with five colours.\\
($D_{2}$) There are some maps for which three colours are insufficient.\\
($D_{3}$) A computer has analysed every type of planar map and verified that each of them is 4-colorable.
\end{quotation}
of the claim $C$, that
\begin{quotation}
\noindent ($C$) Four colours suffice to colour any planar map.
\end{quotation}
by employment of the warrant $W$, which has backing $B$
\begin{quotation}
\noindent ($W$) The computer has been properly programmed and its hardware has no defects.\\
($B$) Technology and computer programming are sufficiently reliable. \cite[pp.~142f.]{Alcolea}
\end{quotation}
He regards this as making clear that, since the warrant is not wholly mathematical, the proof must leave open the possibility of `a specific counterexample, that is to say, a particular map that cannot be coloured with four colours might still exist' \cite[p.~143]{Alcolea}.\footnote{`\dots un contraexemple espec\'ific, \'es a dir, que es trobe un mapa particular que no puga colarar-se amb quatre colors'}

This example demonstrates both the strengths and the dangers of this approach.  To complete Toulmin's layout we are obliged to make explicit not merely the premises and the conclusion, but also the nature of the support which the former is supposed to lend the latter.  Thus the the focus of Appel and Haken's critics, the heterodox deployment of a computer in a mathematical proof, is made glaringly obvious.  However, it is premature to draw from this surface dissimilarity the inference that Appel and Haken's result is less convincing than other mathematical proofs.  A closer reading of Alcolea's reconstruction may clarify this point.  Premises $D_{1}$ and $D_{2}$ have conventional mathematical proofs, as Alcolea points out.  ($D_{1}$ is not strictly relevant to the derivation of $C$, although its proof originated techniques which were instrumental to Appel and Haken's work.)  $D_{3}$ is a very concise summary of the central results of Appel and Haken's work.  It may help to spell out the details at greater length.

There are two essential ideas behind the Appel and Haken proof: unavoidability and reducibility.  An unavoidable set is a set of configurations, that is countries or groups of adjacent countries, at least one of which must be present in any planar map.  For example, all such maps must contain either a two-sided, a three-sided, a four-sided or a five-sided country, so these configurations constitute an unavoidable set.  A configuration is reducible if any map containing it may be shown to be four-colorable.  Two-sided, three-sided, and four-sided countries are all reducible.  To prove the four colour theorem it suffices to exhibit an unavoidable set of reducible configurations.  Alfred Kempe, who introduced the concepts of unavoidability and reducibility, was believed to have proved the four colour theorem in 1879 by showing that five-sided countries were also reducible \cite{Kempe}.  However, in 1890 a flaw was discovered in his reasoning: the five-sided country is not reducible, hence a larger unavoidable set is required if all its configurations are to be reducible.  Appel and Haken used a computer to search for such a set, eventually discovering one with 1,482 members.  The unavoidability of this set could be demonstrated by hand, but the reducibility of all its members would be far too protracted a task for human verification.  Subsequent independent searches have turned up other unavoidable sets.  The smallest to date is a set of 633 reducible configurations found by Neil Robertson, Daniel Sanders, Paul Seymour and Robin Thomas in 1994.\footnote{For details of their publications, see \cite[p.~244]{Wilson}.}  Verifying the reducibility of these configurations still requires a computer.

So for there to be a non-four-colourable planar map, as Alcolea suggests, Appel and Haken (and their successors) must have erred either in the identification of the unavoidable set, or in the demonstration of the reducibility of its member configurations.  Since the former step can be verified by conventional methods, the computer can only be suspected of error in demonstrating reducibility.   Two sorts of computer error should be distinguished: a mistake may be made in the programming, or a fault may arise in the computer itself (the hardware or firmware).  The former error would arise due to a human failure to correctly represent the mathematical algorithms which the computer was programmed to implement.  This sort of mistake does not seem to be interestingly different from the traditional type of mathematical mistake, such as that made by Kempe in his attempt to prove the four colour conjecture.  The second sort of error is genuinely new.  However, it would seem to be profoundly unlikely.

Computer hardware can exhibit persistent faults, some of which can be hard to detect.\footnote{For example, the notorious Pentium bug.}.  However, the potential risks of such faults can be minimized by running the program on many different machines.  One might still worry about Appel and Haken's programs, since they were written in machine code and would therefore be implemented in more or less exactly the same manner on any computer capable of running them, perhaps falling foul of the same bug each time.  This sort of checking might be suspected of being no better than buying two copies of the same newspaper to check the veracity of its reporting.\footnote{As Wittgenstein once remarked in a different context: \cite[\S 265]{PI}.}  However, the same reducibility results were achieved independently, using different programs, as part of the refereeing process for Appel and Haken's work.  Moreover, the more recent programs of Robertson \textit{\& al.}\ were written in higher level languages, as are the programs employed in most other computer-assisted proofs.  The existence of different compilers and different computer platforms ensures that these programs can be implemented in many intrinsically different ways, reducing the likelihood of hardware or firmware induced error to the astronomical.

Thus we may derive an alternative reconstruction of Appel and Haken's argument:
``Given that ($D_{4}$) the elements of the set $U$ are reducible, we can ($Q$) almost certainly claim that ($C$) four colours suffice to colour any planar map, since ($W$) $U$ is an unavoidable set (on account of ($B$) conventional mathematical techniques), unless ($R$) there has been an error in either (i) our mathematical reasoning, or (ii) the hardware or firmware of all the computers on which the algorithm establishing $D_{4}$ has been run.''
If, in addition, we observe that (i) appears to be orders of magnitude more likely than (ii), then $C$ would seem to be in much less doubt than it did in the light of Alcolea's reconstruction.  The purpose of the preceding has been not so much to rescue the four colour conjecture from Alcolea's critique (although few if any graph theorists would accept that a counterexample is possible), but to show up the limitations of Toulmin's pattern as a descriptive technique.  As other critics have pointed out, reconstructing an argument along Toulmin's lines `forces us to rip propositions out of context' \cite[p.~318]{Willard}.  The degree of abstraction necessary to use the diagram at all can make different, incompatible, reconstructions possible, leaving the suspicion that any such reconstruction may involve considerable (and unquantified) distortion.

\section{Walton's new dialectic}
There has been significant progress in informal logic since the publication of \textit{The uses of argument}.  One milestone was the publication of Charles Hamblin's \textit{Fallacies} \cite{Hamblin} in 1970.   This demonstrated the inadequacies of much of traditional fallacy theory and, by way of remedy, proposed an influential dialectical model of argumentation.  Further impetus has come from the recent work of communication theorists such as Frans van Eemeren and Rob Grootendorst \cite{vanEemeren++}.  One contemporary logician who shows the influence of both traditions is Douglas Walton.\footnote{Walton has published a great number of works on informal logic.   \cite{w98} provides an overview of the general method common to many of them.}  The focus of his work is the dialectical context of argument.  
Walton distinguishes between `inference', defined as a set of propositions, one of which is warranted by the others, `reasoning', defined as a chain of inferences, and `argument', defined as a dialogue employing reasoning.  This dialectical component entails that arguments require more than one arguer: at the very least there must be an assumed audience, capable in principle of answering back.  

Winston Churchill once praised the argumentational skills of the celebrated barrister and politician F.~E.~Smith, 1st Earl of Birkenhead, by stressing their suitability to context: `The bludgeon for the platform; the rapier for a personal dispute; the entangling net and unexpected trident for the Courts of Law; and a jug of clear spring water for an anxious perplexed conclave' \cite[p.~176]{Churchill}.  Toulmin also stresses the domain specificity of good practice in argument.  What is distinctive about Walton's analysis is the attempt to characterize dialectical context in terms of general features which are not themselves domain specific.  Without pretending to have an exhaustive classification of argumentational dialogue, he is able to use these features to draw several important distinctions.  The principal features with which he is concerned are the `initial situation' and the `main goal' of the dialogue.  The initial situation describes the circumstances which give rise to the dialogue, in particular the differing commitments of the interlocutors.  The main goal is the collective outcome sought by both (all) participants, which may be distinct from their individual goals.

\begin{table}[h]
\caption{Walton \& Krabbe's `Systematic survey of dialogue types' \protect\cite[p.~80]{wk}}
\label{survey}
\begin{center}
\begin{footnotesize}
\begin{tabular}{ll}
\raisebox{-.55in}{\rotatebox{90}{Main Goal}}&
\begin{tabular}{|p{1.4in}|p{0.75in}|p{0.75in}|p{0.75in}|}
\multicolumn{1}{c}{}&\multicolumn{3}{c}{Initial Situation} \smallskip\\
\cline{2-4}
\multicolumn{1}{l|}{}&Conflict&Open Problem&Unsatisfactory \mbox{Spread of} \mbox{Information}\\
\hline
\raggedright Stable Agreement/ Resolution&Persuasion&Inquiry&Information Seeking\\
\hline
\raggedright Practical Settlement/ Decision (Not) to Act&Negotiation&Deliberation&\multicolumn{1}{c}{}\\
\cline{1-3}
\raggedright Reaching a (Provisional) Accommodation&Eristic&\multicolumn{2}{c}{}\\
\cline{1-2}
\end{tabular}
\end{tabular}
\end{footnotesize}
\end{center}
\end{table}

If we simplify the situation by permitting each discussant to regard some crucial proposition as either true, false or unknown, four possibilities emerge.  Either (0) the discussants agree that the proposition is true (or that it is false), in which case there is no dispute; or (1) one of them takes it to be true and the other false, in which case they will be in direct conflict with each other; or (2) they both regard it as unknown, which may result in a dialogue as they attempt to find out whether it is true or false; or (3) one of them believes the proposition to be true (or false) but the other does not know which it is.
Thus we may distinguish three types of initial situation from which an argumentational dialogue may arise: a conflict, an open problem, or an unsatisfactory spread of information.
A conflict may produce several different types of dialogue depending on how complete a resolution is sought.  For a stable outcome one interlocutor must persuade the other, but, even if such persuasion is impossible they may still seek to negotiate a practical compromise on which future action could be based.  Or they may aim merely to clear the air by expressing their contrasting opinions, without hoping to do more than merely agree to disagree: a quarrel.  These three goals---stable resolution, practical settlement and provisional accommodation---can also be applied to the other two initial situations, although not all three will be exemplified in each case.  So open problems can lead to stable resolutions, or if this is not achievable, to practical settlement.  However, provisional accommodation should not be necessary if the problem is genuinely open, since neither discussant will be committed to any specific view.  Where the dialogue arises merely from the ignorance of one party then a stable resolution should always be achievable, obviating the other goals.  The interplay of these different types of initial situation and main goal thus allows Walton to identify six principal types of dialogue, Persuasion, Negotiation, Eristic, Inquiry, Deliberation and Information Seeking, which may be represented diagrammatically as in Table \ref{survey}.  The contrasting properties of these different types of dialogue are set out in Table \ref{walton}.  This table also states the individual goals of the interlocutors typical to each type, and includes two derivative types: the debate, a mixture of persuasion and eristic dialogue, and the pedagogical dialogue, a subtype of the information seeking dialogue.
Many other familiar argumentational contexts may be represented in terms of Walton's six basic types of dialogue by such hybridization and subdivision.\footnote{See Table 3.1 in \cite[p.~66]{wk} for some further examples.}  

\begin{table}[t]
\caption{Walton's types of dialogue \protect\cite[p.~605]{w97}}
\label{walton}
\begin{footnotesize}
\begin{center}
\begin{tabular}{||p{0.75in}||p{0.75in}|p{0.75in}|p{0.75in}|p{0.75in}||}
\hline
\hline
\textbf{\mbox{Type~of} \mbox{Dialogue}}&\textbf{Initial \mbox{Situation}}&\textbf{Individual \mbox{Goals~of} \mbox{Participants}}&\textbf{Collective \mbox{Goal~of} \mbox{Dialogue}}&\textbf{Benefits}\\
\hline
\hline
Persuasion&Difference of opinion&Persuade other party&Resolve difference of opinion&Understand positions\\
\hline
Inquiry&Ignorance&Contribute findings&\raggedright Prove or disprove conjecture&Obtain \mbox{knowledge}\\
\hline
Deliberation&Contemplation of future consequences&Promote personal goals&Act on a thoughtful basis&Formulate personal priorities\\
\hline
Negotiation&Conflict of interest&Maximize gains (self-interest)&Settlement (without undue inequity)&Harmony\\
\hline
Information-Seeking&One party lacks information&Obtain information&Transfer of knowledge&Help in goal activity\\
\hline
\raggedright Quarrel (Eristic)&\raggedright Personal conflict&Verbally hit out at and humiliate opponent&Reveal deeper conflict&Vent emotions\\
\hline
Debate&Adversarial&Persuade third party&Air strongest arguments for both sides&Spread \mbox{information}\\
\hline
Pedagogical&Ignorance of one party&Teaching and learning&Transfer of knowledge&Reserve \mbox{transfer}\\
\hline
\hline
\end{tabular}
\end{center}
\end{footnotesize}
\end{table}

It is central to Walton's work that the legitimacy of an argument should be assessed in the context of its use: what is appropriate in a quarrel may be inappropriate in an inquiry, and so forth.  Although some forms of argument are never legitimate (or never illegitimate), most are appropriate if and only if they are ``in the right place''.  For example, threats are inappropriate as a form of persuasion, but they can be essential in negotiation.  In an impressive sequence of books, Walton has analyzed a wide variety of fallacious or otherwise illicit argumentation as the deployment of strategies which are sometimes admissible in contexts in which they are inadmissible.  However, Walton has not directly addressed mathematical argumentation.  In the next section I shall set out to explore how well his system may be adapted to this purpose.

\section{Applying Walton to mathematics}

In what context (or contexts) do mathematical proofs occur?  The obvious answer is that mathematical proof is a special case of inquiry.  Indeed, Walton states that the collective goal of inquiry is to `prove or disprove [a] conjecture'.  An inquiry dialogue proceeds from an open problem to a stable agreement.  That is to say from an initial situation of mutual ignorance, or at least lack of commitment for or against the proposition at issue, to a main goal of shared endorsement or rejection of the proposition.  This reflects a standard way of reading mathematical proofs: the prover begins from a position of open-mindedness towards the conjecture, shared with his audience.  He then derives the conjecture from results upon which they both agree, by methods which they both accept.

But this is not the only sort of dialogue in which a mathematical proof may be set out.  As William Thurston has remarked, mathematicians `prove things in a certain context and address them to a certain audience' \cite[p.~175]{Thurston}.  Indeed, crucially, there are several different audiences for any mathematical proof, with different goals.  Satisfying the goals of one audience need not satisfy those of the others.  For example, a proof may be read by:
\begin{itemize}
\item Journal referees, who have a professional obligation to play devil's advocate;
\item Professional mathematicians in the same field, who may be expected to quickly identify the new idea(s) that the proof contains, grasping them with only a few cues, but who may already have a strong commitment to the falsehood of the conjecture;
\item Professional mathematicians in other (presumably neighbouring) fields, who will need more careful and protracted exposition;
\item Students and prospective future researchers in the field, who could be put off by too technical an appearance, or by the impression that all the important results have been achieved;
\item Posterity, or in more mercenary terms, funding bodies: proof priority can be instrumental in establishing cudos with both.
\end{itemize}
This list suggests that the initial situation of a proof dialogue cannot always be characterized as mutual open-mindedness.  Firstly, in some cases, the relationship between the prover and his audience will be one of conflict.  If the conjecture is a controversial one, its prover will have to convince those who are committed to an incompatible view.  And if an article is refereed thoroughly, the referees will be obliged to adopt an adversarial attitude, irrespective of their private views.  

Secondly, as the later items indicate, proofs have a pedagogic purpose.
Thurston relates his contrasting experiences in two fields to which he made substantial contributions.  As a young mathematician, he proved many results in foliation theory using powerful new methods.  However, his proofs were of a highly technical nature and did little to explain to the audience how they too might exploit the new techniques.  As a result, the field evacuated: other mathematicians were afraid that, by the time they had mastered Thurston's methods, he would have proved all the important results.  In later work, on Haken manifolds, he adopted a different approach.  By concentrating on proving results which provided an infrastructure for the field, in a fashion which allowed others to acquire his methods, he was able to develop a community of mathematicians who could pursue the field further than he could alone.  The price for this altruism was that he could not take all the credit for the major results.
Proofs which succeed in the context Thurston advocates proceed from an initial situation closer to Walton's `unsatisfactory spread of information'.  This implies that information seeking is another context in which mathematical proofs may be articulated.  Of course, the information which is being sought is not merely the conjecture being proved, but also the methods used to prove it.  

An unsatisfactory spread of information, unlike a conflict or an open problem, is an intrinsically asymmetrical situation.  We have seen that proofs can arise in dialogues wherein the prover possesses information sought by his interlocutors.  Might there be circumstances in which we should describe as a proof a dialogue in which the prover is the information seeker?
This is the question considered by Thomas Tymoczko \cite[p.~71]{Tymoczko} and Yehuda Rav \cite{Rav}, to somewhat different ends.  Tymoczko considers a community of Martian mathematicians who have amongst their number an unparalleled mathematical genius, Simon.  Simon proves many important results, but states others without proof.  Such is his prestige, that ``Simon says'' becomes accepted as a form of proof amongst the Martians.  Rav considers a fantastical machine, Pythiagora, capable of answering mathematical questions instantaneously and infallibly.  Both thought experiments consider the admission of a dialogue with an inscrutable but far better informed interlocutor as a possible method of proof.  

In both cases we are invited to reject this admission, although, interestingly, for different reasons.  Rav sees his scenario as suggesting that proof cannot be purely epistemic: if it were then Pythiagora would give us all that we needed, but Rav suggests that we would continue to seek conventional proofs for their other explanatory merits.  He concludes that Pythiagora could not give us proof.  Tymoczko draws an analogy between his thought experiment and the use of computers in proofs such as that of the four colour theorem.  He argues that there is no formal difference between claims backed by computer and claims backed by Simon: they are both appeals to authority.  The difference is that the computer can be a warranted authority.  Hence, on Tymoczko's admittedly controversial reading, computer assisted proof is an information seeking dialogue between the prover and the computer.

So far we have seen that the initial situation of a proof dialogue can vary from that of an inquiry.  What of the main goal---must this be restricted to stable resolution?  Some recent commentators have felt the need for a less rigorous form of mathematics, with a goal closer to Walton's practical settlement.  Arthur Jaffe and Frank Quinn \cite{Jaffe+} introduced the much discussed, if confusingly named, concept of `theoretical mathematics'.  They envisage a division of labour, analogous to that between theoretical and experimental physics, between conjectural or speculative mathematics and rigorous mathematics.  Where traditional, rigorous mathematicians have theorems and proofs, theoretical mathematicians make do with `conjectures' and `supporting arguments'.  This echoes an earlier suggestion by Edward Swart \cite{Swart} that we should refrain from accepting as theorems results which depend upon lengthy arguments, whether by hand or computer, of which we cannot yet be wholly certain.  He suggests that `these additional entities could be called agnograms, meaning theoremlike statements that we have verified as best we can but whose truth is not known with the kind of assurance that we attach to theorems and about which we must thus remain, to some extent, agnostic' \cite[p.~705]{Swart}.  In both cases the hope is that further progress will make good the shortfall: neither Jaffe and Quinn's conjectures nor Swart's agnograms are intended as replacements for rigorously proved theorems.

More radical critics of the accepted standards of mathematical rigour suggest that practical settlement can be a goal of proof and not merely of lesser, analogous activities.  For instance, Doron Zeilberger \cite{Zeilberger} envisages a future of semi-rigorous (and ultimately non-rigorous) mathematics in which the ready availability by computer of near certainty reduces the pursuit of absolute certainty to a low resource allocation priority.  Hence he predicts that a mathematical abstract of the future could read ``We show, in a certain precise sense, that the Goldbach conjecture is true with probability larger than $0.99999$, and that its complete truth could be determined with a budget of $\$ 10$ billion'' \cite[p.~980]{Zeilberger}.  `Proofs' of this sort explicitly eschew stable resolution for practical settlement.  Thus Zeilberger is arguing that proofs could take the form of deliberation or negotiation.

The last of Walton's dialogue types is the eristic dialogue, in which no settlement is sought, merely a provisional accommodation in which the commitments of the parties are made explicit.  This cannot be any sort of proof, since no conclusion is arrived at.  But it is not completely without interest.  A familiar diplomatic euphemism for a quarrel is ``a full and frank exchange of views'', and such activity does have genuine merit.  Similarly, even failed mathematical proofs can be of use, especially if they clarify previously imprecise concepts, as we saw with Kempe's attempted proof of the four colour conjecture \cite{Kempe}.  This process has something in common, if not with a quarrel, at least with a debate, which we saw to be a related type of dialogue.

\begin{table}[h]
\caption{Some types of proof dialogue}
\label{proof}
\begin{footnotesize}
\begin{center}
\begin{tabular}{||p{0.75in}||p{0.75in}|p{0.75in}|p{0.75in}|p{0.75in}||}
\hline
\hline
\textbf{\mbox{Type~of} \mbox{Dialogue}}&\textbf{Initial \mbox{Situation}}&\textbf{Main Goal}&\textbf{\mbox{Goal~of} \mbox{Prover}}&\textbf{\mbox{Goal~of} \mbox{Interlocutor}}\\
\hline
\hline
\raggedright Proof as Inquiry&Open-mindedness&Prove or disprove conjecture& \raggedright Contribute to outcome& Obtain \mbox{knowledge}\\
\hline
\raggedright Proof as Persuasion& \raggedright Difference of opinion& Resolve difference of opinion with rigour& \raggedright Persuade interlocutor&Persuade prover\\
\hline
\raggedright Proof as Information-Seeking (Pedagogical)&\raggedright Interlocutor lacks information&\raggedright Transfer of knowledge&\raggedright Disseminate knowledge of results \& methods&Obtain \mbox{knowledge}\\
\hline
 \raggedright `Proof' as Information-Seeking (\textit{e.g.} Tymoczko)& \raggedright Prover lacks information& \raggedright Transfer of knowledge&Obtain \mbox{information}&Presumably \mbox{inscrutable}\\
\hline
 \raggedright `Proof' as Deliberation (\textit{e.g.} Swart)&Open-mindedness& \raggedright Reach a provisional conclusion& \raggedright Contribute to outcome&Obtain \mbox{warranted} belief\\
\hline
 \raggedright `Proof' as Negotiation (\textit{e.g.} Zeilberger)& \raggedright Difference of opinion& \raggedright Exchange resources for a provisional conclusion&Contribute to outcome&Maximize value \mbox{of exchange}\\
\hline
 \raggedright `Proof' as Eristic/ Debate&Irreconcilable difference of opinion&Reveal deeper conflict&Clarify \mbox{position}&Clarify \mbox{position}\\
\hline
\hline
\end{tabular}
\end{center}
\end{footnotesize}
\end{table}

To take stock, we have seen that most of Walton's dialogue types are reflected to some degree in mathematical proof.  Table \ref{proof}, an adaptation of Table \ref{walton}, sets out the difference between the various types of proof dialogue introduced.

\section{Proof dialogues}
In this last section I shall explore how the classification of proof dialogues may help to clarify many of the problems that have arisen in the  philosophical debate over the nature of mathematical proof.  We can see that proofs may occur in several distinct types of dialogue, even if we do not count the suspect cases (the entries for Table \ref{proof} where `proof' is in scare quotes).  An ideal proof will succeed within inquiry, persuasion and pedagogic proof dialogues.  Suboptimal proofs may fail to achieve the goals of at least one of these dialogue types.  In some cases, this may be an acceptable, perhaps inevitable, shortcoming; in others it would fatally compromise the argument's claim to be accepted as a proof.

As Thurston's experience with foliation theory demonstrated, not every proof succeeds pedagogically.  Proofs in newly explored areas are often hard to follow, and there are some results which have notoriously resisted all attempts at clarification or simplification.\footnote{For example, von Staudt's proof of the equivalence of analytic and synthetic projective geometry has retained its difficulty for nearly two centuries.  See \cite[pp.~193 f.]{Rota} for a discussion.}  Yet, if these proofs succeed in inquiry and persuasion dialogues, we have no hesitation in accepting them.  Conversely, there are some `proofs' which have a heuristic usefulness in education, but which would not convince a more seasoned audience.  Pedagogic success is neither necessary nor sufficient for proof status---but it is a desirable property, nonetheless.

An argument might convince a neutral audience, but fail to persuade a determined sceptic.  Just this happened to Andrew Wiles's first attempt at a proof of the Fermat conjecture: the initial audience were convinced, but the argument ran into trouble when exposed to determined criticism from its referees.  Such a case might be seen as success within an inquiry proof dialogue, followed by failure in a persuasion proof dialogue.  A similar story could be told about Kempe's `proof' of the four colour conjecture: a result which received far less scrutiny than Wiles's work, and was thereby widely accepted for eleven years.  On the other hand, if even the sceptics are convinced, then an open-minded audience should follow suit.  Thus, on the conventional understanding of mathematical rigour, success within both inquiry and persuasion proof dialogues is necessary for an argument to count as a proof.

We saw in the last section how a variety of differently motivated departures from the prevailing standards of mathematical rigour may be understood as different types of proof dialogue.  One of Walton's principal concerns in his analysis of natural argumentation is the identification of shifts from one type of dialogue to another.  Such shifts can take a variety of forms: either gradual or abrupt, and either replacing the former type of dialogue or embedding the new type within the old.  These processes are an essential and productive aspect of argumentation, but they are also open to abuse.  Similar warnings apply to the less rigorous types of proof dialogue.

Many of the concerns which the critics of these forms of argumentation have advanced may be understood as an anxiety about illicit shifts of proof dialogue type.  For example, the published discussion of mathematical conjecture is something which Jaffe and Quinn welcome: their concern is that such material not be mistaken for theorem-proving.  Although some of their critics interpreted their advocacy of `theoretical mathematics' as a radical move, their primary goal was a conservative one: to maintain a sharp demarcation between rigorous and speculative work.  Their `measures to ensure ``truth in advertising''\thinspace' \cite[p.~10]{Jaffe+} are precisely calculated to prevent illicit shifts between inquiry and deliberation proof dialogues.\footnote{Indeed, misleading advertisements are one of Walton's principal examples of an illicit dialogue shift.}  A similar story could be told about Tymoczko or Swart's discussion of methods they see as falling short of conventional rigour.  Zeilberger is advocating the abandonment of rigour, but he recognizes at least a temporary imperative to separate rigorous from `semi-rigorous' mathematics.

As Toulmin \textit{\& al.}\ remark `it has never been customary for philosophers to pay much attention to the \textit{rhetoric} of mathematical debate' \cite[p.~89]{Toulmin+}.  The goal of this article has been to exhibit some of the benefits that may accrue from a similarly uncustomary interest in informal logic.

\end{document}